\documentclass[12pt]{article}
\makeatletter

\@addtoreset{equation}{section}
\makeatother
\usepackage{latexsym}
\makeatletter
\AtBeginDocument{\@ifpackageloaded{natbib}{\ifNAT@numbers\if@filesw\immediate\write\@auxout{\string\global\string\NAT@numberstrue}\fi\fi}{}}
\makeatother
\begin{document}
\renewcommand{\thefootnote}{\fnsymbol{footnote}}
\renewcommand{\theequation}{%
\thesection.\arabic{equation}}

  \title{ Correction
    to \\
     ``The classification of the surfaces with parallel mean curvature vector in two-dimensional complex space 
     forms"}
\author{Katsuei Kenmotsu}
 \date{}

\setcounter{footnote}{0}
\newtheorem{lm}{Lemma}
\newtheorem{corr}{Correction}
\newtheorem{prop}{Proposition}
\newtheorem{thm}{Theorem}
\newtheorem{cor}{Corollary}
\newtheorem{remark}{Remark}
\maketitle
Abstract: We give a condition under which the findings of the paper cited above work well and determine
  the surfaces that  were not considered before.

\section{Correction}
Kenmotsu and Zhou \cite{kenzhou} claimed  that the  surfaces stated in the
title were locally classified. The proof
used the system of ordinary differential equations  in Ogata \cite{ogata}, but in that paper \cite{ogata},
 a mistake was made, as pointed out by Hirakawa \cite{hirakawa}. In fact,
the claim that ``$\lambda$ is a real-valued function defined on $U$" made in line 3 of  page 401
in \cite{ogata}  is not correct in general. Hence,  our derivation \cite{kenzhou}  needs an additional
assumption: {\it $\lambda = \bar{\lambda}$}.  We follow the notation used in \cite{kenzhou} and \cite{ogata}.
 In  \cite{hirakawa}, Hirakawa classified the surfaces with $a=\bar{a}$,
 which is weaker than  {\it $\lambda = \bar{\lambda}$}. Hence, we were forced to study the case 
 of $a \neq \bar{a}$ for the classification of these surfaces.\footnote[0]{Partly supported by 
 JSPS Grant-in-Aid for Scientific 
 Research (C-25400062)} 
  In this paper, we prove the following:
\vspace{0.3cm}

{\bf Correction} {\it A  parallel mean curvature vector
surface  in a complex two-dimensional   complex space
form with  $a \neq \bar{a}$   depends on   one real-valued
harmonic function on the surface  and five real constants  if the ambient space is not flat, the
mean curvature vector does not vanish, and the Kaehler angle  is not constant.}
\vspace{0.3cm}

We remark that  Hirakawa \cite{hirakawa} proved the above result  under an additional condition $c \equiv 0$ (for $c$ see \S $2$).
As a by-product of the correction, we have the following corollary:
\vspace{0.3cm}

{\bf Corollary}
{\it Any two-dimensional smooth manifold can be locally embedded in the complex projective plane 
and in the
complex hyperbolic plane as a parallel mean curvature vector surface.}
\vspace{0.3cm}

  Some comments about the assumptions  of the correction are in order. Chern and Wolfson \cite{cherwolf} and
  Eschenburg, Guadalupe, and Tribuzy \cite{escgt}  studied
  the case of a zero  mean curvature vector, i.e.,  minimal surfaces.
 The case of a constant Kaehler angle is included in  Hirakawa \cite{hirakawa}, because the
 Gaussian curvature of such a surface is constant.   Parallel mean curvature vector surfaces
 in a flat complex space
  form are
   locally classified independently  by Hoffman \cite{hoffman}, Yau \cite{yau}, and
   Chen \cite{chen1}, because
   the ambient space is locally isometric to $R^{4}$ with  the standard metric.

\section{Analysis of the structure equations}
To  correct our calculation, we rewrite the structure equations of the immersion.
  Let $\overline{M}[4\rho]$ be a complex two-dimensional complex space form with a constant holomorphic
  sectional curvature $4\rho$,
$M$ be an oriented and connected real two-dimensional Riemannian manifold with  Gaussian curvature
$K,$ and $x: M\longrightarrow
\overline{M}[4\rho]$ be an isometric immersion with  Kaehler angle  $\alpha$ such that the  mean
curvature vector field $H$  is nonzero and parallel for the normal connection on the normal bundle
of the
immersion. In this section, we show that the Kaehler angle plays an important role in understanding
the immersion.

Since all calculations and formulas  in \cite{ogata} are valid until page 400,  there exists a local
field of unitary coframes
$\{w_{1},w_{2}\}$ on $\overline{M}[4\rho]$ such that, by restricting it to $x$, the Riemannian metric $ds^2$
 on $M$ is written as  $ds^{2}=\phi\bar{\phi}$, where
$\phi=\cos\alpha/2\cdot\omega_{1}+\sin\alpha/2\cdot\bar{\omega}_{2}$.  Let $a$ and $c$ be the complexified second fundamental tensors of $x$ with respect to $\{ \omega_{1}, \omega_{2}\}$.
 Then, the Kaehler angle $\alpha$ and the complex 1-form $\phi$  satisfy
\begin{eqnarray}
d\alpha &=& ( a +  b )\phi + (\bar{a} + b )\bar{\phi} , \\
d\phi &=&  (\bar{a} - b)\cot \alpha \cdot \phi \wedge \bar{\phi},
\end{eqnarray}
where $2b = |H| >0$. By (2.4), (2.5), and (2.6) of \cite{ogata}, the Gauss, Codazzi-Mainardi, and Ricci equations of $x$ are, respectively,
 \begin{eqnarray}
&& K= -4(|a|^{2} - b^{2})+6\rho\cos^{2}\alpha, \\
&& da\wedge\phi= -\left(2a (\bar{a} - b) \cot\alpha  + \frac{3}{2}\rho\sin \alpha \cos \alpha \right)
\phi\wedge\bar{\phi}, \\
&& dc\wedge\bar{\phi}=2c(a - b)\cot\alpha \cdot \phi\wedge\bar{\phi}, \\
&& |c|^{2} = |a|^{2} + \frac{\rho}{2}(-2 + 3\sin^{2}\alpha).
\end{eqnarray}
Conversely, for a real number  $\rho$, a positive number $b$, a complex-valued 1-form $\phi$, a
real-valued function
$\alpha$, complex-valued functions $a$ and $c$ on a  simply connected domain $D$ in $R^2$,
if they  satisfy (2.1)--(2.6), where  $K$ denotes the Gaussian curvature of the
 Riemannian metric
$ds^2 =\phi \bar{\phi}$ on $D$, then there exists
an isometric immersion $x : (D, ds^{2}) \longrightarrow \overline{M}[4\rho]$ such that  the
Kaehler angle of $x$
is  $\alpha$, the mean curvature
vector of $x$ is parallel, and its length is equal to  $2b$.

The immersion $x$ is called  a general type if it satisfies $ d\alpha \neq 0$, $a \neq \bar{a}$,
and $c \neq 0$ on $M$.
From now on, we suppose that the ambient space is nonflat and  study  the immersions
of  a general type.
We remark that the results in this paper are local in nature, and we suppose that all formulas
hold in
a neighborhood of $M$.

We  define $a_{1}$ and $c_{1}$  by
\begin{equation}
\left. \bigg\{
\begin{array}{@{\,}ll}
da = (a_{1} + at_{1})  \phi + t_{2}\bar{\phi},  \\
dc = 2c(a- b)\cot\alpha\cdot \phi + (c_{1} + c\overline{t_{1}}) \bar{\phi},
\end{array}
\right.
\end{equation}
where $t_{i}= t_{i}(\alpha,a,\bar{a})\; (i=1,2,\dots)$ are  functions of $\alpha$, $a$, and $\bar{a}$
listed in the Appendix.

Taking the exterior differentiation of $(2.6)$, we have
 \begin{equation}
    c\overline{c_{1}} = \bar{a}a_{1}.
\end{equation}
The exterior differentiation of $(2.7)$ allows us to define  $a_{11}$ and $c_{11}$ by
\begin{equation}
\bigg\{
\begin{array}{@{\,}ll}
 da_{1} = a_{11} \phi + (3a_{1}(\bar{a}-b) \cot \alpha  + t_{3}) \bar{\phi},  \\
 dc_{1} = (3c_{1}(a-b)\cot \alpha  + c t_{4})\phi + c_{11}\bar{\phi}.
 \end{array}
\end{equation}
The exterior differentiation of $(2.8)$  implies, by (2.7) and (2.9),
\begin{equation}
|c_{1}|^{2} = |a_{1}|^{2} + t_{5}.
\end{equation}
By (2.8), we have $|c|^2 |c_{1}|^2 = |a|^2 |a_{1}|^2$. This with (2.6) and (2.10) implies
$$
\frac{\rho}{2}(-2 + 3\sin^{2}\alpha)|a_{1}|^2 = -(|a|^2 + \frac{\rho}{2}(-2 + 3\sin^{2}\alpha)t_{5}.
$$
Hence, if $\rho \neq 0$ and $\sin^{2}\alpha \neq 2/3$, then
\begin{equation}
|a_{1}|^{2} = t_{6}.
\end{equation}
We note that $a_{1} \neq 0$, since $x$ is of a general type. Indeed, if $a_{1}$ vanishes identically, then so do
$c_{1}$, and hence,  $t_{4}$.  However,  by the definition of $t_{4}$, we calculate
$ t_{4} - \overline{t_{4}} = -6b(a-\bar{a})(5+ 3\cos 2\alpha)\sin^2 \alpha/(1+3 \cos 2\alpha)^2$,
 which gives a contradiction.

Taking the exterior differentiation of (2.11),  we have
\begin{equation}
 a_{11} = \frac{t_{7} a_{1} + t_{8}}{\overline{a_{1}}}.
\end{equation}
Inserting this into the first equation of (2.9), we take the exterior differentiation of the formula, which
 gives us
\begin{equation}
t_{9} a_{1} + \overline{t_{9} a_{1}} + t_{10} =0.
\end{equation}
We study first the case of  $t_{9} \neq 0$. By (2.11) and (2.13), $a_{1}$ is determined
in terms of $\alpha$, $a$, and $\bar{a}$, denoted
$a_{1}=t_{11}(\alpha,a,\bar{a})$.
Let us take the exterior differentiation of this formula. Then, by (2.9) and (2.12) we have  two
equations
  for $\alpha$, $a$, and $\bar{a}$, denoted
  $t_{12}(\alpha,a, \bar{a})=0$ and $t_{13}(\alpha,a, \bar{a})=0$.   When $t_{9}(\alpha,a, \bar{a})=0$,
  we have $t_{10}(\alpha,a, \bar{a})=0$ by (2.13).
  These are all relations for $\alpha$, $a$, and $\bar{a}$ that we can  get
  by  exterior differentiations of the structure equations.
We show the following lemma:
\begin{lm} If $x$ is of a general type, then
 $d\alpha \wedge da =0$.
 \end{lm}
 Proof. First, we prove Lemma 1 when $t_{9} \neq 0$.  It holds that
 $\partial t_{12}/\partial a \neq 0$ or $\partial t_{12}/\partial \bar{a} \neq 0$. Indeed,
 if both are zero, then $t_{12}$ is a function of $\alpha$ and then $t_{12}=0$ implies that
 $\alpha$ is constant, giving a contradiction. We suppose that $\partial t_{12}/\partial \bar{a} \neq 0$.
 Then, $t_{12}(\alpha,a,\bar{a})=0$ can be solved as $\bar{a}=\bar{a}(\alpha,a)$ by the implicit
 function theorem. Inserting this into the expression for $t_{13}$ gives $t_{13}(\alpha,a,\bar{a}(\alpha,a))=0$.
 This proves Lemma 1. When $\partial t_{12}/\partial a \neq 0$, we prove Lemma 1 in the same
 way as before, because $\alpha$ is real valued.
In the case of $t_{9} = 0$ also, we can prove Lemma 1  in the
same manner as before. Indeed, $t_{12}$ and $t_{13}$  now replace
$t_{9}$ and $t_{10}$.
This completes the proof of  Lemma 1.




Last, we study $c$ in terms of  $\alpha$, $a$, and $\phi$.
  From (2.7) and  (2.8),
   we have
  \begin{equation}
  d \log c = 2(a-b) \cot \alpha \cdot \phi + \left(\frac{c_{1}}{c} + \overline{t_{1}}\right)\bar{\phi},
  \end{equation}
where
\begin{equation}
\frac{c_{1}}{c} = \frac{\bar{c}c_{1}}{\bar{c}c}
   = \frac{a \overline{t_{11}}(\alpha,a,\bar{a})}{|a|^2 + \rho/2(-2 + 3 \sin^2\alpha)}.
   \end{equation}
  Since
    its length $|c|$ is uniquely determined by (2.6),
   we have shown the following:
   \begin{lm} $c$ is determined up
   to  a real constant by integration of the formula $(2.14)$.
   \end{lm} 

\section{Proof of the correction}
  If $\alpha$ is not constant, then
  $a$ is a function of $\alpha$, as given by $a=a(\alpha)$, from Lemma 1.
  It follows
 from (2.1) and (2.7) that
 \begin{equation}
 \frac{da}{d\alpha} = \frac{\cot \alpha}{\overline{a(\alpha)} +b}\left( - 2ba(\alpha) +  2|a(\alpha)|^{2}
 + \frac{3\rho}{2}\sin^{2}\alpha \right), \ (a+b \neq 0).
 \end{equation}
 The complex-valued function $a=a(\alpha)$ is determined by the above first-order ordinary differential
 equation; hence, it depends on one complex constant.
 Next, we determine  the Kaehler angle  $\alpha$ and the Riemannian metric $ds^{2}$ on $M$.  Let
$\phi = \lambda(z,\bar{z})dz,$ where $z$ is an isothermal coordinate on $M$. From (2.1), we have
$\alpha_{z}=\lambda (a+b)$; hence $a + b \neq 0$.  Taking its exterior derivative, (2.7) implies
$\alpha_{z\bar{z}} =
\lambda_{\bar{z}}(a+b) + \lambda\bar{\lambda}t_{2}$. From (2.2), we see that
 $\lambda_{\bar{z}} = - (\bar{a}-b) \lambda\bar{\lambda}\cot \alpha$. From these formulas,
 we have
 \begin{equation}
\alpha_{z\bar{z}} - F(\alpha)\alpha_{z}\alpha_{\bar{z}} =0,
\end{equation}
where
$$
F(\alpha) = \frac{((a(\alpha) - b)(\overline{a(\alpha)} - b)
+ 3\rho/2 \sin^{2} \alpha)}{(a(\alpha) + b)(\overline{a(\alpha)} + b)} \cot \alpha .
$$

\begin{lm}
 Any solution $\alpha$ of  $(3.2)$ is  written 
as $\alpha(z,\bar{z}) = \psi (f(z,\bar{z}))$,
where
$f(z,\bar{z})$ is  a real-valued harmonic function on $M$ and $\psi$ is  a real solution of the second-order
  ordinary differential equation 
\begin{equation}
\psi''(t) - F(\psi)\psi'(t)^{2} =0.
\end{equation}
\end{lm}

 Proof. 
Define a real valued function $K(t)$ of one real variable  by
$$
K(t) = \int e^{-\int F(t)dt}dt
$$
and set $f(z,\bar{z})=K(\alpha(z,\bar{z}))$. By $(3.2)$, $f(z,\bar{z})$ is a harmonic function, i.e., $f$ satisfies $\partial^{2}f/\partial z\partial \bar{z} =0$.
We set  $\psi(t) = K^{-1}(t)$. Then, $\psi(t)$  satisfies $(3.3)$, proving Lemma 3. 

Because $\psi(t)$ depends on  two real constants,  the Kaehler angle $\alpha$ is determined by one harmonic function on $M$
and two real constants.
By (2.1),   the Riemannian metric $ds^2$ is written as
$$
ds^2= \left|\frac{\alpha_{z}}{a(\alpha)+b} \right|^{2}dzd\bar{z}.
$$
In view of Lemmas 1 and  2, we proved the following:
 \begin{thm}  Let $x: M\longrightarrow
\overline{M}[4\rho]\ (\rho \neq 0)$  be a  parallel mean curvature vector  immersion such that
$H \neq 0$. If $x$ is
of a general type, then it is determined by one
real-valued harmonic function on $M$
 and five real constants up to isometries of $\overline{M}[4\rho]$.
 \end{thm}
 Now we discuss  the converse of the above and  prove the existence of the immersion of a general type.
  Given  $\rho \neq 0$ and $b>0$,  we take a solution $a=a(\alpha)$ of (3.1). For
  any real-valued harmonic function $f(z,\bar{z})$ on a simply connected domain $D \subset R^{2}$
  with $f_{z} \neq 0$ and a solution $\psi$ of (3.3), we define
  $\alpha = \psi(f(z,\bar{z}))$.  We set
  $\lambda = \alpha_{z}/(a(\alpha) +b), \
  \phi =\lambda dz$, and $ds^2= \phi \bar{\phi}$.  It holds that $d\phi  = (\bar{a} - b )
  \cot \alpha \cdot \phi \wedge \bar{\phi}$.
These $\alpha$, $a=a(\alpha)$, $\phi$, and $ds^2$ satisfy $(2.1)$--$(2.4)$.  $a_{1}$
 is defined by (2.7); hence, we have
 $$
 a_{1}= -a(\alpha) t_{1}(\alpha, a(\alpha),\overline{a(\alpha)}) + (a(\alpha) + b)\frac{da}{d\alpha}.
 $$
We determine  $c$ explicitly as follows:
We set
$$
\omega = \frac{1}{2i}\cdot \frac{\omega_{1} \phi - \overline{\omega_{1}\phi}}
{|a(\alpha)|^2 + \rho/2(-2 + 3 \sin^2\alpha)},
$$
where
\begin{eqnarray*}
\omega_{1} &=& \left(|a(\alpha)|^2
+ \frac{\rho}{2}(-2 + 3\sin^2\alpha) \right)\left(2(a(\alpha)-b)\cot \alpha
- \overline{t_{1}}(\alpha, a(\alpha), \overline{a(\alpha)}) \right)  \nonumber \\
  &&    + \overline{a(\alpha)}\left(a(\alpha) t_{1}(\alpha, a(\alpha),\overline{a(\alpha)})
  - (a(\alpha) + b)\frac{da}{d\alpha} \right).
\end{eqnarray*}
Then, we have the following:
\begin{lm}
$\omega$ is a closed real 1-form on $D$.
\end{lm}
We shall omit the proof, because
this is verified by direct calculation using (2.1), (2.2), (2.4), and (3.1).

According to Lemma 4, there exists a real-valued function $\nu =\nu(z,\bar{z})$ on
$D$ such that $d\nu = \omega$. Set
\begin{equation}
c= (|a(\alpha)|^2 + \frac{\rho}{2}(-2 + 3 \sin^{2}\alpha)^{1/2}e^{i\nu(z,\bar{z})}.
\end{equation}
Since this satisfies (2.5) and (2.6), these $\alpha$, $a$, $\phi$, and $c$ satisfy the structure
equations $(2.1)$--$(2.6)$. We
 have thus  proved the following:
\begin{thm}   Given  $\rho \neq 0, \ b>0$, and a nonconstant real-valued
 harmonic function $f$
 on a simply connected domain $D$ in $R^{2}$,
there exist a  Riemannian metric $ds^2$ on $D$ and
 a parallel mean curvature vector immersion from $D$ into $\overline{M}[4\rho]$ of a general type
  such that  $|H|=2b$, and the Kaehler angle is determined by $f$.
\end{thm}
The correction is proved by Theorems $1$ and  $2$.
\vspace{0.5cm}

{\bf Proof of  Corollary}\quad  Let $(u,v)$ be a  local coordinate of a two-dimensional smooth manifold $M$.
Then, the coordinates  are  nonconstant harmonic functions on the neighborhood of $M$.
Hence, we  prove the corollary by  setting  $f=u$ or $f=v$ in Theorem 2.

\vspace{0.5cm}
\section{Associated family}
We show the explicit example of  the second fundamental tensors  of the surfaces  with
 $a \neq \bar{a}$ and $c \neq 0$.
  For real numbers
$c_{1},c_{2}$, we set
\begin{eqnarray*}
a(t) &=& \frac{-4+(9+4c_{1})\sin^{2}t - 9c_{1}\sin^{4}t
+ \sqrt{-1} \sqrt{2(8-9\sin^{2}t)(-1+ c_{1}\sin^{2}t)}}
{4(-1+ c_{1}\sin^{2}t) - \sqrt{-1}\sqrt{2(8-9\sin^{2}t)(-1+ c_{1}\sin^{2}t)}},  \\
\xi(t) &=& 2^{5/2} \int \frac{\cot t}{\sqrt{(8-9\sin^2 t)(-1 +c_{1}\sin^2 t)}}dt + c_{2}, \\
c(t) &=& \sqrt{\frac{c_{1}}{2(-9+8c_{1})}} (8 -9 \sin^{2}t)e^{\sqrt{-1} \xi(t)}.
\end{eqnarray*}
Given   a nonconstant real-valued harmonic function $f(z,\bar{z})$
on
a domain $D$ in $R^2$, we set $\alpha = f$ and $\phi = f_{z}/(a(f)+1) dz$. Then, for any $c_{1}$
satisfying $c_{1} <0 $ or $c_{1} > 9/8$,
and any $c_{2} \in R$,
these $\alpha$,  $a(f)$, $c(f)$, and $\phi$ satisfy the
structure equations $(2.1)$--$(2.6)$
with $b=1$, and $\rho=-3$. Hence, they define a two-parameter family of
parallel mean curvature vector
surfaces
in $\overline{M}[-12]$ with $a \neq \bar{a}$ and $c \neq 0$ such that the Kaehler angle
is  $f$. In particular, if we change the value of $c_{2}$ under a fixed $c_{1}$, then
the resulting surfaces are isometric and have the same length for the mean curvature vectors.
Hence, these give us the associated family of parallel mean curvature vector surfaces
in $\overline{M}[-12]$.

We note that the case of $c_{1}=0$ in the example above is already found in Hirakawa \cite{hirakawa}.

\section{Appendix} 
{\footnotesize 
\begin{eqnarray*}
&& t_{1}(\alpha,a,\bar{a}) = ( -4b + 
 12b \sin^2 \alpha +4a +3a \sin^2\alpha ) \frac{\cot \alpha}{-2 + 3 \sin^2 \alpha},  \\
&& t_{2}(\alpha,a,\bar{a}) =   2a(\bar{a}-b) \cot \alpha  
+ \frac{3}{2}\rho \sin \alpha \cos \alpha ,  \\
&& t_{3}(\alpha,a,\bar{a}) = -t_{1}t_{2} - t_{2} (a-b)\cot\alpha + 3a(\bar{a}-b)t_{1} \cot \alpha 
 - a(\bar{a} +b)\frac{\partial t_{1}}{\partial \alpha} - at_{2}\frac{\partial t_{1}}{\partial a} \\
&& \qquad \qquad \qquad+(a+b) \frac{\partial t_{2}}{\partial \alpha} 
 + \overline{t_{2}}\frac{\partial t_{2}}{\partial \bar{a}},  \\
 && t_{4}(\alpha,a,\bar{a}) = 2\left(t_{2} \cot \alpha - (a-b)(\bar{a}-b) \cot^2\alpha 
 - \frac{(a-b)(\bar{a}+b)}{\sin^2\alpha} \right) + \overline{t_{1}}(a-b) \cot \alpha \\ 
 && \qquad \qquad \qquad  - (a+b)\frac{\partial\overline{t_{1}}}{\partial \alpha}  
  - \overline{t_{2}}\frac{\partial \overline{t_{1}}}{\partial \bar{a}}, \\  
&& t_{5}(\alpha,a,\bar{a}) = t_{3}\bar{a} -\overline{t_{4}} \left(a\bar{a} 
+ \frac{\rho}{2}(-2+3\sin^2\alpha)\right),   \\ 
&& t_{6}(\alpha, a,\bar{a}) = \frac{-t_{5}}{\rho/2 (-2 + 3 \sin^{2} \alpha) } ( a \bar{a} 
+ \frac{\rho}{2} (-2 + 3 \sin^{2} \alpha) ), \\
&& t_{7}(\alpha,a,\bar{a}) = -\overline{t_{3}} + \frac{\partial t_{6}}{\partial a},  \\
&& t_{8}(\alpha,a,\bar{a}) = - 3 t_{6} (a - b) \cot \alpha + (a + b) \frac{\partial t_{6}}{\partial \alpha} 
+ a t_{1} \frac{\partial t_{6}}{\partial a} + \overline{t_{2}} \frac{\partial t_{6}}{\partial \bar{a}},  \\
&& t_{9}(\alpha,a,\bar{a}) = t_{6} \left(-(\bar{a} - b) t_{7} \cot \alpha 
+ (\bar{a} + b)  \frac{\partial t_{7}}{\partial \alpha} 
+ t_{2} \frac{\partial t_{7}}{\partial a} + \overline{a t_{1}} \frac{\partial t_{6}}{\partial \bar{a}} \right) 
- t_{7} \overline{t_{8}},  \\
&& t_{10}(\alpha,a,\bar{a}) = t_{6} \left( t_{3} t_{7}  - t_{7} \overline{t_{7}}
-4 (\bar{a} - b) t_{8} \cot \alpha 
+ (\bar{a} + b) \frac{\partial t_{8}}{\partial \alpha} + t_{2} \frac{\partial t_{8}}{\partial a} 
+ \overline{a t_{1}}\frac{\partial t_{8}}{\partial \bar{a}} \right) } 
- t_{8} \overline{t_{8} \\
&& \qquad  \qquad  \qquad - t_{6} ^2 \left(3 \overline{t_{2}}  \cot \alpha
- 3 \frac{(\bar{a} - b) (a + b)}{\sin^{2} \alpha} - 3 (\bar{a} - b) (a - b) \cot^{2} \alpha 
+ \frac{\partial t_{3}}{\partial a} - \frac{\partial t_{7}}{\partial \bar{a}} \right), \\
&& t_{11}(\alpha,a,\bar{a}) =  \frac{1}{2t_{9}} \left( -t_{10} \pm \left(t_{10}^2 
- 4 |t_{9}|^{2}t_{6} \right)^{1/2} \right), \\
&& t_{12}(\alpha,a,\bar{a}) = t_{7}t_{11} + t_{8} - t_{6}\frac{\partial t_{11}}{\partial a} 
- \overline{t_{11}} \left( (a + b) \frac{\partial t_{11}}{\partial \alpha} 
+ a t_{1} \frac{\partial t_{11}}{\partial a} 
+ \overline{t_{2}} \frac{\partial t_{11}}{\partial \bar{a}}   \right), \\
&& t_{13}(\alpha,a,\bar{a}) =t_{3} + 3 (\bar{a}-b)t_{11} \cot \alpha 
-\left( (\bar{a} + b) \frac{\partial t_{11}}{\partial \alpha} 
+ t_{2} \frac{\partial t_{11}}{\partial a} 
+ \overline{t_{11}} \frac{\partial t_{11}}{\partial \bar{a}}  
 + \overline{at_{1}}\frac{\partial t_{11}}{\partial \bar{a}} \right), \\
&&  \bar{t_{i}}(\alpha,a, \bar{a}) = t_{i}(\alpha,\bar{a},a) \quad (i=1,2,3,...) . 
\end{eqnarray*}

\medskip
\begin{flushleft}

 Katsuei Kenmotsu \\
Mathematical Institute,  Tohoku University  \\
980-8578 \quad  Sendai, Japan \\
email:  kenmotsu@math.tohoku.ac.jp
\end{flushleft}
\end{document}